\documentclass[12pt]{article}
\usepackage[english]{babel}
\usepackage{Custom}
\usepackage{AuxiliaryFiles/AuxiliaryFiles}
\usepackage{Custom-2}
\begin{document}
\title{\vspace{-.8cm} \Title\vspace{-.65cm}}
\author{\Author\thanks{Chair of Algebra and Number Theory, RWTH Aachen University, Pontdriesch 14--16, 52062 Aachen, Germany. E-mail: \EMail}}
\date{\vspace{-.2cm} \Date}
\maketitle
\UnderlineBox{\Abstract}

\medskip
\noindent\textbf{Mathematics Subject Classification (MSC 2020):}
14N10 (primary); 12H05, 26C10, 14Q20, 13P15 (secondary).

\vspace{-.25cm}
\section*{Introduction}
Let \(k\) be a field of characteristic zero, and fix an algebraic closure \(K\) of \(k\).
Consider \(k\) as a differential field endowed with the trivial derivation \(\partial\) and work in the corresponding differential polynomial ring \(k\{x_1,x_2\}\).
We are interested in the differential system
\EqL{\Sigma_2\colon \begin{cases}
	x_1^\prime-x_2^2=0\\
	x_2^\prime-x_1^2=0
\end{cases}}{SystemStudied}
The natural object associated with \(\Sigma_2\) is the differential ideal
\[I=[x_1^\prime-x_2^2, x_2^\prime-x_1^2]\]
It is the smallest ideal of \(k\{x_1,x_2\}\) containing the polynomials defining the equations of \(\Sigma_2\) and that is closed under \(\partial\): for example we have
\Eq{(\smash{x_1^\prime-x_2^2})^\prime=x_1^{\prime\prime}-2x_2x_2^\prime\in I\\
(\smash{x_2^\prime-x_1^2})^\prime=x_2^{\prime\prime}-2x_1x_1^\prime\in I}
or, in other words, the equations
\EqL{x_1^{\prime\prime}-2x_2x_2^\prime=0\\
x_2^{\prime\prime}-2x_1x_1^\prime=0}{FirstDerivative}
are differential consequences of the equations in \(\Sigma_2\).
In the algebraic setting, these new equations need not follow from the original ones, because the jet symbols are treated as algebraically independent variables.
The algebraic system
\[\Sigma_2^{(1)}\colon \begin{cases}
	x_1^\prime-x_2^2=0\\
	x_2^\prime-x_1^2=0\\
	x_1^{\prime\prime}-2x_2x_2^\prime=0\\
	x_2^{\prime\prime}-2x_1x_1^\prime=0
\end{cases}
\]
obtained by adjoining the two equations \eqref{FirstDerivative} to \(\Sigma_2\) involves six variables, so its solution set is an algebraic variety in \(K^6\).

\paragraph*{}
More generally, let \(n\in\N\), \(n\ge2\) be given. Consider the system
\EqL{\Sigma_{n}\colon\begin{cases} x_1^\prime-x_2^n=0 \\ x_2^\prime-x_1^n=0 \end{cases}}{SigmaN}
and write \(f_1=x_1^\prime-x_2^n\) and \(f_2=x_2^\prime-x_1^n\) for ease of notation.
For a fixed \(m\in\N_0\) we shall denote by \(\Sigma_{n}^{(m)}=\Prol{m}{\Sigma_{n}}\) the algebraic system consisting of all the equations of the form
\[\partial^hf_l=0 \qquad\text{for some }h\in\Set{0,\dots,m} \text{ and }l\in\Set{1,2}\]
Thus \(\Sigma_{n}^{(m)}\) consists of the equations of \(\Sigma_{n}\) together with their derivatives up to order \(m\).
We denote the corresponding algebraic variety by
\[\X_{n}^{(m)}=Z(\smash{\Sigma_{n}^{(m)}})\subseteq\A_{K}^{2(m+2)}\]

\paragraph*{}
All computational experiments reported in this paper are carried out over \(k=\Q\) using OSCAR \cite{oscar}.
A Jupyter notebook containing the construction of the prolongation systems, executable examples, and instructions for reproducing the computations is available in the GitHub repository \cite{MyGitHub2}.

\paragraph*{}
The second-order equations in \(\Sigma_2^{(1)}\) can be simplified using the two equations in \(\Sigma_2\).
We obtain the equivalent system\vspace{-.12cm}
\Eq{\begin{cases}x_1^\prime-x_2^2=0\\
x_2^\prime-x_1^2=0\\
x_1^{\prime\prime}-2x_1^2x_2=0\\
x_2^{\prime\prime}-2x_1x_2^2=0
\end{cases}}
Generalizing this procedure to \(\Sigma_n^{(m)}\), we can express every derivative of \(x_1\) and \(x_2\) as a polynomial in \(x_1\) and \(x_2\) using the lower order equations of the system.
For every \(n\ge2\) and \(m\ge0\), the algebraic variety \(\X_{n}^{(m)}\) is equivalently defined by
\EqL{\begin{cases}
x_1^{(l)}&=A_{n,l}(x_1,x_2)\\
x_2^{(l)}&=B_{n,l}(x_1,x_2)\quad\text{for all }l\in\Set{1,\dots,m+1}
\end{cases}
}{EqP}
for suitable polynomials \(A_{n,l},B_{n,l}\in k[x_1,x_2]\).
To emphasize this property, we write \(x_1=a\) and \(x_2=b\) when considering \(x_1\) and \(x_2\) as parameters for the algebraic variety \(\X_{n}^{(m)}\).

\paragraph*{}
By the chain rule, the equations of \(\Sigma_n\) induce the following derivation on \(k[a,b]\):
\[
D_n=b^n\frac{\partial}{\partial a}
  +a^n\frac{\partial}{\partial b}
\]
The polynomials \(A_{n,l}\) and \(B_{n,l}\) are determined recursively by the initial values \(A_{n,0}=a\), \(B_{n,0}=b\) and the relations
\Eq{A_{n,l+1}&=D_n(A_{n,l})\\
B_{n,l+1}&=D_n(B_{n,l})\qquad\quad\text{for all }l\in\N_0}
By induction on \(l\), it is immediate to check that \(A_{n,l}\) and \(B_{n,l}\) are homogeneous polynomials of degree \(1+l(n-1)\) with integer coefficients.
\subparagraph*{}
This discussion leads to the following parametrization of \(\X_{n}^{(m)}\).
For simplicity, since the field \(K\) is fixed, we write \(\A^N\) and \(\P^N\) for the affine and projective spaces over \(K\) of dimension \(N\) when no ambiguity arises.

\begin{Pro}{}{Parametrization}
The algebraic variety \(\X_{n}^{(m)}\subseteq\A^{2(m+2)}\) is the graph of the polynomial map
\[
\NMapH{\A^2}{\A^{2(m+1)}}{(a,b)}{(A_{n,1}(a,b),B_{n,1}(a,b),\dots,A_{n,m+1}(a,b),B_{n,m+1}(a,b))}{0pt}
\]
In particular
\[
\MapH{\Phi_{n}^{(m)}}{\A^2}{\X_n^{(m)}}{(a,b)}{(a,b,A_{n,1}(a,b),B_{n,1}(a,b),\dots,A_{n,m+1}(a,b),B_{n,m+1}(a,b))}{0pt}
\]
is an isomorphism of algebraic varieties, and \(\X_n^{(m)}\) is smooth, irreducible of dimension 2 for every \(n\ge2\) and \(m\ge0\).
\end{Pro}

\begin{Proof}
The first statement follows directly from the above discussion.
The inverse of \(\Phi_n^{(m)}\) is the restriction to \(\X_n^{(m)}\) of the projection onto the first two coordinates.
\end{Proof}

\paragraph*{}
Recall that the \textbf{geometric degree} \(\deg\X\) of an irreducible affine variety \(\X\subseteq\A^N\) of dimension \(d\) is the number of points in its intersection with a generic affine linear subspace of the ambient space \(\A^N\) of codimension \(d\).
Since \(\X_{n}^{(m)}\) has dimension 2 by the preceding proposition, its geometric degree can be computed by intersecting it with two generic affine hyperplanes of \(\A^{2(m+2)}\).
Although \(\X_{n}^{(m)}\simeq\A^2\) for every \(n\ge2\), the geometric degree of \(\X_{n}^{(m)}\) depends on its embedding into \(\A^{2(m+2)}\) and may vary with \(n\) and \(m\).

\paragraph*{}
The problem of computing geometric degrees of prolongation varieties is largely open. In his recent work \cite[Proposition 7.7]{Lanciano2026}, Leonardo Lanciano obtained bounds for these degrees when the original system satisfies certain regularity conditions.
For a polynomial dynamical system, Gleb Pogudin asked whether the geometric degree of its prolongation varieties is always a quasi-polynomial function of the prolongation order \cite[Section 2, Question 7]{Pogudin2025}. In particular, he proposed the following formula for the system \(\Sigma_2\):

\begin{Conj}[before upper={\setlength{\abovedisplayskip}{0pt}\setlength{\abovedisplayshortskip}{-12pt}}]{}{OriginalConjecture}
\begin{equation}
\deg\smash{\X_2^{(m)}}=
\begin{cases}
(m+2)^2-2 & \text{if }m\equiv 1\bmod{3}\\
(m+2)^2   & \text{otherwise}
\end{cases}\label{Conjecture}
\end{equation}
\end{Conj}

\noindent
The aim of this paper is to prove this conjecture.
The ideas behind the proof originated from computational experiments in OSCAR using algebraic Thomas decomposition.
The next section explains how these experiments led to a criterion for computing the geometric degree through the analysis of the intersections at infinity.
The proof is then given in Section \ref{ProofConjecture}.
In Section \ref{Thomas}, we present algebraic Thomas decomposition and derive a formula for counting the distinct points in finite intersections of a prolongation variety with pairs of affine hyperplanes.
Finally, Section \ref{Further} discusses extensions to the family \(\Sigma_n\) and presents a general conjecture for the geometric degree of the corresponding prolongation varieties.
\begin{Not*}{}{}
In Sections \ref{Motivation} and \ref{ProofConjecture}, we only consider the quadratic system \(\Sigma_2\).
For ease of notation, we will write \(\X^{(m)}\) and \(\Phi^{(m)}\) in place of \(\X_{\smash{2}}^{(m)}\) and \(\Phi_{\smash{2}}^{(m)}\). We shall also write \(D\), \(A_l\) and \(B_l\) in place of \(D_2\), \(A_{2,l}\) and \(B_{2,l}\).
\end{Not*}

\section{Computational motivation}\label{Motivation}
The geometric idea behind the proof of Conjecture \ref{conj:OriginalConjecture} came from experiments in OSCAR, where algebraic Thomas decomposition was used to count the number of distinct intersection points of a prolongation variety with pairs of affine hyperplanes.
To determine whether this quantity coincided with the geometric degree, we studied the intersections at infinity of the plane curves obtained by substituting \(\Phi^{(m)}\) into the hyperplane equations.
\subparagraph*{}
Fix \(m\in\N_0\) and consider two affine hyperplanes of \(\A^{2(m+2)}\):
\[
H_1\colon c_0+\sum_{h=0}^{m+1} \left(c_h^A x_1^{(h)}+c_h^B x_2^{(h)}\right)=0\qquad\qquad
H_2\colon d_0+\sum_{h=0}^{m+1}
\left(d_h^A x_1^{(h)}+d_h^B x_2^{(h)}\right)=0
\]
When computing the geometric degree, we may restrict to hyperplanes with coefficients in \(k\), since the \(k\)-rational points are Zariski-dense in the parameter space of hyperplanes.
By Proposition~\ref{pro:Parametrization}, the intersection \(\X^{(m)}\cap H_1\cap H_2\) corresponds to the solution set in \(\A^2\) of the polynomial system
\[
\begin{cases}
\displaystyle P(a,b)=c_0+\sum_{h=0}^{m+1} \left(c_h^A A_h+c_h^B B_h\right)=0\\[12pt]
\displaystyle Q(a,b)=d_0+\sum_{h=0}^{m+1} \left(d_h^A A_h+d_h^B B_h\right)=0
\end{cases}
\]
For a generic choice of the hyperplanes, this system has \(\deg\X^{(m)}\) distinct solutions in \(K^2\).
Since we may restrict to a smaller nonempty Zariski-open subset without losing genericity, we additionally require that the matrix
\[
\begin{pmatrix}
c_{m+1}^A & c_{m+1}^B\\
d_{m+1}^A & d_{m+1}^B
\end{pmatrix}\in k^{2\times 2}
\vspace{.12cm}\]
be invertible.
We can then replace \(P\) and \(Q\) by suitable linear combinations without changing their common zero set, so that
\[P=A_{m+1}+U \qquad\text{and}\qquad Q=B_{m+1}+V\]
where\vspace{-.25cm}
\[U,V \in W_m = \text{Span}_k(1,A_0,B_0,\dots,A_m,B_m)\]
are polynomials of degree at most \(m+1\).
In our OSCAR experiments, we took \(U,V\) as linear combinations of \(1,A_0,B_0,\dots,A_m,B_m\) with random integer coefficients. For pairs \(P,Q\) with finitely many common zeros, we counted these points using the method presented in Section \ref{Thomas}.

\paragraph*{}
The point count need not equal the geometric degree, since the chosen hyperplanes may not be generic.
For example, when \(m=1\), the conjectured geometric degree is \(\deg\X^{(1)}=7\).
Consider the two affine hyperplanes with equations
\Eq{
H_1\colon\quad&
x_1^{\prime\prime}+5x_2^\prime+2x_1^\prime
+5x_2-4x_1+5=0\\
H_2\colon\quad&
x_2^{\prime\prime}-x_1^\prime-4x_2-2x_1+4=0
}
Setting \(x_1=a\) and \(x_2=b\) as in the introduction, we have the system
\Eq{\begin{cases}P(a,b)=2a^2b+5a^2+2b^2+5b-4a+5=0\\
Q(a,b)=2ab^2-b^2-4b-2a+4=0\end{cases}}
in \(\Q[a,b]\).
Using the method described in Section \ref{Thomas} and the counting formula \eqref{FormulaCountingPoints} we obtain six distinct affine solutions.
To explain why this choice gives six affine points rather than seven, work in the projective plane \(\P^2_K\) and homogenize \(P,Q\) introducing a new coordinate \(Z\):
\Eq{P^h(Z,a,b)&=2a^2b+(\smash{5a^2+2b^2})Z+(5b-4a)Z^2+5Z^3\\
Q^h(Z,a,b)&=2ab^2-b^2Z-(4b+2a)Z^2+4Z^3}
The restrictions of \(P^h\) and \(Q^h\) to \(Z=0\) are the homogeneous parts of degree 3 of \(P\) and \(Q\), respectively \(2a^2b\) and \(2ab^2\).
Thus, the projective curves have exactly two common points at infinity, \([0:1:0]\) and \([0:0:1]\).
Modifying the lower-degree terms of \(P,Q\) leaves the points at infinity unchanged, although their intersection multiplicities may change.
\subparagraph*{}
In particular, we can add a term \(\lambda a^2\) to \(Q\) with \(\lambda\in\Q\).
Let us then consider the polynomial\vspace{-.35cm}
\[Q_\lambda=Q+\lambda a^2\]
In the chart \(a\ne 0\),
\Eq{P^h(Z,1,b)&=5Z+2b+\text{(terms of larger degree)}\\
Q_\lambda^h(Z,1,b)&=\lambda Z+\text{(terms of larger degree)}
}
The Jacobian determinant at \(p=[0:1:0]\) is
\[\det{\begin{pmatrix}
	\partial_Z P^h & \partial_b P^h \\
	\partial_Z Q_\lambda^h & \partial_b Q_\lambda^h \\
\end{pmatrix}}_{(0,0)} = \det\begin{pmatrix}
	5 & 2 \\ \lambda & 0
\end{pmatrix}=-2\lambda\]
We conclude that the intersection at \(p\) is transversal if and only if \(\lambda\ne0\).
A similar computation shows that the intersection at \(q=[0:0:1]\) is transversal for every \(\lambda\).
Transversality at these two fixed points is expressed by the non-vanishing of Jacobian determinants and is therefore an open condition on the hyperplane coefficients. Since it holds for \(\lambda\ne0\), we conclude that both points at infinity have intersection multiplicity 1 for a generic choice of the hyperplanes.
\noindent
Replacing \(Q\) by \(Q_1=Q+a^2\) gives a system with seven distinct affine solutions, as verified using algebraic Thomas decomposition. Together with the above analysis at infinity, this proves the conjecture for \(m=1\).

\subsection{Geometric degree and intersections at infinity}\label{Infinity}
We extend this analysis to arbitrary prolongation orders and derive a criterion for computing the geometric degree.
Let \(m\) be a non-negative integer, and write \(l=m+1\) and \(d=m+2\) for ease of notation.
As above, we may assume that the restrictions of two generic affine hyperplanes to \(\X^{(m)}\) have the form
\[
P=A_l+U\qquad\text{and}\qquad Q=B_l+V
\]
for some \(U,V\in W_m\).
The polynomials \(A_l\) and \(B_l\) are homogeneous of degree \(d\), whereas every element of \(W_m\) has degree at most \(d-1\), so both \(P\) and \(Q\) have degree \(d\).
Write
\[
P=F+F_{d-1}+\cdots+F_0\qquad\text{and}\qquad
Q=G+G_{d-1}+\cdots+G_0
\]
where \(F_h\) and \(G_h\) are homogeneous of degree \(h\), whilst \(F=A_l\) and \(G=B_l\).
Introducing a new variable \(Z\), the homogenizations of \(P\) and \(Q\) are
\[
P^h=F+ZF_{d-1}+\cdots+Z^dF_0
\qquad\text{and}\qquad
Q^h=G+ZG_{d-1}+\cdots+Z^dG_0
\]
On the line at infinity \(Z=0\), their common zeros satisfy
\EqL{
F(a,b)=G(a,b)=0
}{CommonZero}
Consequently, the points at infinity are determined by the homogeneous polynomial
\[
H=\gcd(F,G)=\gcd(A_l,B_l)
\]
In particular, \(H\) does not depend on the coefficients of \(U\) and \(V\).
Suppose that \(H\) is square-free, and let \(r=\deg H\).
Then \(H\) has exactly \(r\) distinct zeros in \(\P_K^1\) and the two projective curves have \(r\) common points at infinity.

\paragraph*{}
Following the case of the problematic hyperplanes, we are interested in determining when the intersections at these points are transversal.
Consider the homogeneous polynomials\vspace{-.2cm}
\[
J_a=\frac{\partial G}{\partial a}F_{d-1}-\frac{\partial F}{\partial a}G_{d-1}
\qquad\text{and}\qquad
J_b=\frac{\partial G}{\partial b}F_{d-1}-\frac{\partial F}{\partial b}G_{d-1}
\vspace{.08cm}
\]
In the affine chart \(b\ne0\), setting \(b=1\), the Jacobian determinant with respect to the local coordinates \((Z,a)\) of the homogenized equations at \(Z=0\) is\vspace{.12cm}
\[
\det
\begin{pmatrix}
\partial_ZP^h & \partial_aP^h\\
\partial_ZQ^h & \partial_aQ^h
\end{pmatrix}_{Z=0}
=
J_a(a,1)
\]
Similarly, in the chart \(a\ne0\), the Jacobian determinant with respect to \((Z,b)\) at \(Z=0\) is \(J_b(1,b)\).
Applying Euler's identity to the homogeneous polynomials \(F\) and \(G\) gives
\[
a\frac{\partial F}{\partial a}
+b\frac{\partial F}{\partial b}=dF
\qquad\text{and}\qquad
a\frac{\partial G}{\partial a}
+b\frac{\partial G}{\partial b}=dG
\vspace{.12cm}\]
Using the definition of \(J_a\) and \(J_b\) we obtain
\[
aJ_a+bJ_b=d(GF_{d-1}-FG_{d-1})
\]
Consequently, using Equation \eqref{CommonZero}, we have
\[
aJ_a(a,b)+bJ_b(a,b)=0
\]
at any common point \(p=[0:a:b]\) at infinity.
If \(b\ne 0\), the Jacobian criterion and this identity show that the intersection at \(p\) is non-transversal exactly when \(J_a(a,b)=J_b(a,b)=0\).
If \(b=0\), then \(a\ne0\), and the same criterion follows by exchanging the roles of \(a\) and \(b\).
\subparagraph*{}
Since the two charts cover the line at infinity, we conclude that the non-transversal intersections at infinity are the common projective zeros of \(H,J_a\) and \(J_b\).
These are homogeneous polynomials in two variables, so they have a common projective zero if and only if they have a nonconstant common factor.
All intersections at infinity are transversal if and only if\vspace{-.35cm}
\[\gcd(H,J_a,J_b)=1\]
The intersection multiplicity at each of these points is 1.

\begin{The}{}{TheConj}
Suppose that \(H\) is square-free and that there exists a choice of \(U,V\in W_m\) for which the projective curves of equations \(P^h=0\) and \(Q^h=0\) intersect transversally at every common point at infinity.
Then
\Eq{
\deg\X^{(m)}
&=d^2-r\\[-2pt]
&=(m+2)^2-\deg H
}
\end{The}

\begin{Proof}
Let \(\CC=W_m\times W_m\) be the affine space over \(k\) parametrizing the pairs \((U,V)\), and recall that \(P=A_l+U\), \(Q=B_l+V\).
Since \(H\) is square-free and independent of the coefficients of \(U\) and \(V\), the \(r\) points at infinity are fixed and distinct.
Transversality at these points is an open condition on the coefficients of \(U\) and \(V\): since it holds for at least one choice by assumption, it holds for a generic choice. Thus, the total intersection multiplicity at infinity is generically equal to \(r\).

It remains to show that the affine intersections are generically finite and consist of simple points.
Consider the incidence variety
\[
\II=\{((U,V),(a,b))\in\CC\times\A^2\colon P(a,b)=Q(a,b)=0\}
\]
Since the polynomials \(A_h\) and \(B_h\) are homogeneous of positive degree, we have the direct sum decomposition
\vspace{-.12cm}
\Eq{W_m & =k\cdot 1\oplus W_m^+\\
&=k\cdot 1 \oplus \operatorname{Span}_k(A_0,B_0,\dots,A_{l-1},B_{l-1})}
Writing \(U=c+U_0\) and \(V=e+V_0\) accordingly, the equations defining \(\II\) uniquely determine\vspace{-.4cm}
\Eq{
c=-A_l(a,b)-U_0(a,b)\\
e=-B_l(a,b)-V_0(a,b)
}
Consequently,\vspace{-.24cm}
\[
\II\simeq
W_m^+\times W_m^+\times\A^2
\]
as algebraic varieties.
In particular \(\II\) is smooth and irreducible, with
\[
\dim{}\II
=2\,\dim{}{W_m}
=\dim{}{\CC}
\]

We claim that the natural projection
\[
\pi\colon \II\longrightarrow\CC
\]
is a dominant map.
Indeed, for \(U=a\) and \(V=b\), the origin satisfies \(P=Q=0\).
Since \(A_l,B_l\) are homogeneous of degree \(d\ge2\), their first derivatives vanish at the origin, and therefore
\vspace{-.2cm}
\[
\det
\begin{pmatrix}
\partial_aP & \partial_bP\\
\partial_aQ & \partial_bQ
\end{pmatrix}_{(0,0)}
=1
\]
The restriction of \(\pi\) to a suitable open neighborhood of this point is smooth and hence open \cite[Lemmas 29.35.14 and 29.35.10]{StacksProject}.
The image of \(\pi\) therefore contains a nonempty Zariski-open subset of \(\CC\).
Since \(\CC\) is irreducible, this open subset is dense, and \(\pi\) is dominant.
\\
Since \(\dim{}{\II}=\dim{}{\CC}\), the generic fibers of \(\pi\) are zero-dimensional \cite[Chapter I, Section 6.3]{Shafarevich}.
Moreover, since \(\II\) is smooth and we are working over a field of characteristic zero, the second Bertini theorem \cite[Chapter II, Section 6.2, Theorem 2.27]{Shafarevich} implies that the generic fibers are smooth.
Thus, for a generic choice of \(U\) and \(V\), \(Z(P,Q)\) is finite and every intersection point has multiplicity 1.
\\
We may therefore choose \(U,V\) generically so that the intersections at infinity are transversal, and all affine intersection points are simple.
In particular, the projective curves defined by \(P^h=0\) and \(Q^h=0\) have no common component, and Bézout's theorem gives a total intersection multiplicity of \(d^2\).
Since the \(r\) points at infinity each contribute 1, there are exactly  \(d^2-r\) distinct affine intersection points.
The isomorphism \(\Phi^{(m)}\) identifies \(Z(P,Q)\) with \(\X^{(m)}\cap H_1\cap H_2\): we conclude that
\[
\deg\X^{(m)}=d^2-r
\qedhere
\]
\end{Proof}

\paragraph*{}
Theorem \ref{the:TheConj} determines the geometric degree without requiring an affine point count.
To apply this criterion computationally, consider
\[
P_*=A_l+A_{l-1}\qquad\text{and}\qquad
Q_*=B_l+B_{l-1}
\]
The corresponding polynomials used to check transversality at infinity are
\Eq{
J_a^*=
\frac{\partial B_l}{\partial a}A_{l-1}
-\frac{\partial A_l}{\partial a}B_{l-1}
\qquad\qquad
J_b^*=
\frac{\partial B_l}{\partial b}A_{l-1}
-\frac{\partial A_l}{\partial b}B_{l-1}
}
If \(H=\gcd(A_l,B_l)\) is square-free and \(\gcd(H,J_a^*,J_b^*)=1\), Theorem \ref{the:TheConj} gives
\[
\deg\X^{(m)}=(m+2)^2-\deg H
\]
This criterion involves only polynomial arithmetic and GCD computations.
Using it, we verified the conjectured degree formula in OSCAR for every \(m\in\Set{1,\dots,1000}\).

\paragraph*{}
Comparing Theorem \ref{the:TheConj} with Conjecture \ref{conj:OriginalConjecture} suggests the following statement about \(\gcd(A_l,B_l)\):

\begin{Conj}[before upper={\setlength{\abovedisplayskip}{0pt}\setlength{\abovedisplayshortskip}{-8pt}}]{}{ConjGCD}
\[
\deg\gcd(A_l,B_l)=
\begin{cases}
2&\text{if }l\equiv 2\bmod 3\\
0&\text{otherwise}.
\end{cases}\qquad\quad\text{for all }l\in\N
\]
\end{Conj}

\noindent
In Proposition \ref{pro:EqConj}, we will show that this conjecture implies the hypotheses of Theorem \ref{the:TheConj} and is equivalent to Conjecture \ref{conj:OriginalConjecture}.


\section{The proof of the conjecture}\label{ProofConjecture}
In the introduction, we used the derivation
\[
D=b^2\frac{\partial}{\partial a}+a^2\frac{\partial}{\partial b}
\]
on \(k[a,b]\) to construct the polynomials \(A_l\) and \(B_l\) appearing in the parametrization \(\Phi^{(m)}\). Starting from \(A_0=a\) and \(B_0=b\), these polynomials are defined recursively by
\Eq{A_{l+1}&=D(A_l)\\
B_{l+1}&=D(B_l)\qquad(l\ge 0)}
For a monomial \(a^rb^s\) we have
\[
D(a^rb^s)=ra^{r-1}b^{s+2}+sa^{r+2}b^{s-1}
\]
Note that every nonzero term changes both exponents of \(a\) and \(b\) by \(-1\) modulo 3.
Since \(A_0=a\) and \(B_0=b\), it follows by induction that the exponents \(r,s\) of every monomial in \(A_l\) satisfy\vspace{-.32cm}
\EqL{
r\equiv (1-l)\bmod3
\qquad\text{and}\qquad
s\equiv-l\bmod3
}{CongruenceMod3}
Similarly, the exponents of every monomial in \(B_l\) satisfy
\EqL{
r\equiv-l\bmod3
\qquad\text{and}\qquad
s\equiv(1-l)\bmod3
}{CongruenceMod3-2}
In particular, when \(l\equiv 2\bmod 3\), the exponents of \(a,b\) in each monomial are not congruent to \(0\) modulo 3, and are therefore positive integers. This implies that
\[
ab\mid \gcd(A_l,B_l)
\qquad\quad\text{if }l\equiv 2\bmod3
\]
In this case, \(\deg\gcd(A_l,B_l)=2\) if and only if \(\gcd(A_l,B_l)\) and \(ab\) are associated in \(k[a,b]\).
With \(\sim\) denoting equality up to multiplication by a nonzero scalar, Conjecture \ref{conj:ConjGCD} can be reformulated equivalently as
\EqL{
\boxed{
\gcd(A_l,B_l)\sim
\begin{cases}
ab&\text{if }l\equiv 2\bmod 3\\
1&\text{otherwise}
\end{cases}
\qquad\quad\text{for all }l\in\N }}{ConjFormula}

\begin{Pro}{}{EqConj}
Conjectures \ref{conj:OriginalConjecture} and \ref{conj:ConjGCD} are equivalent.
Moreover, Conjecture \ref{conj:ConjGCD} implies that the hypotheses of Theorem \ref{the:TheConj} hold for every \(m\in\N_0\).
\end{Pro}

\vspace{.12cm}
\begin{Proof}
Fix \(m\in\N_0\) and set \(l=m+1\).
Suppose first that Conjecture \ref{conj:ConjGCD} holds.
If \(l\not\equiv2\bmod3\), then\vspace{-.45cm}
\[
H=\gcd(A_l,B_l)\sim1
\]
Hence there are no common points at infinity, and both hypotheses of Theorem \ref{the:TheConj} are satisfied.

If \(l\equiv2\bmod3\) then \(H\sim ab\), so \(H\) is square-free, and the only common
points at infinity are \(p=[0:1:0]\) and \(q=[0:0:1]\).
Consider the choice
\[
P=A_l+A_{l-1}\qquad\qquad Q=B_l+B_{l-1}
\]
We show that the corresponding projective curves intersect transversally at both points.
The recurrence defining \(A_l\) and \(B_l\) implies that the coefficients of the monomials \(a^lb\) in \(A_l\) and \(a^l\) in \(B_{l-1}\) are nonzero.
Consequently,
\[
\frac{\partial A_l}{\partial b}(1,0)\ne0
\qquad\text{and}\qquad
B_{l-1}(1,0)\ne0
\]
The congruence conditions on the exponents imply that
\[
\frac{\partial B_l}{\partial b}(1,0)=0
\]
Therefore\vspace{-.25cm}
\[
J_b(1,0)=-\frac{\partial A_l}{\partial b}(1,0)B_{l-1}(1,0)
\ne0
\vspace{.12cm}
\]
and the intersection at \(p\) is transversal.
By symmetry, the intersection at \(q\) is also transversal.
Theorem \ref{the:TheConj} therefore applies for \(m\in\N_0\) and gives 
\[
\deg\X^{(m)}
=
(m+2)^2-\deg\gcd(A_{m+1},B_{m+1})
\]
By Conjecture \ref{conj:ConjGCD}, this formula is precisely the one in Conjecture~\ref{conj:OriginalConjecture}.

Conversely, suppose that Conjecture \ref{conj:OriginalConjecture} holds, and write \(d=l+1\).
For a generic choice of the hyperplanes, the plane curves \(P=0\) and \(Q=0\) meet in \(\deg\X^{(m)}\) distinct affine points, all of intersection multiplicity 1.
By Bézout's theorem, the total projective intersection multiplicity of their homogenizations is \(d^2\).
Thus the total intersection multiplicity at infinity is\vspace{-.25cm}
\[
d^2-\deg\X^{(m)}
\]
Let \(H=\gcd(A_l,B_l)\).
Each zero of \(H\) of multiplicity \(e\) contributes at least \(e\) to the total intersection multiplicity at infinity.
Summing over all points at infinity, we obtain
\[
\deg H\le d^2-\deg\X^{(m)}
\]
If \(l\not\equiv2\bmod3\), Conjecture \ref{conj:OriginalConjecture} gives\vspace{-.12cm}
\[
\deg\X^{(m)}=d^2
\]
Thus \(\deg H=0\), and hence \(H\sim1\).

If \(l\equiv2\bmod3\), the congruences proved above imply that \(ab\mid H\), and therefore
\[
\deg H\ge2
\]
On the other hand, Conjecture \ref{conj:OriginalConjecture} gives
\[
\deg\X^{(m)}=d^2-2
\]
and hence \(\deg H\le2\).
It follows that \(\deg H=2\) and \(H\sim ab\).
\end{Proof}

\paragraph*{}
We have reduced the proof of Conjecture \ref{conj:OriginalConjecture} to showing \eqref{ConjFormula}.
We can perform a further simplification, reformulating this GCD formula as a single coprimality condition on suitable univariate polynomials.
Work in the chart \(a\ne0\) and introduce the coordinate \(z=b/a\). Since \(A_l,B_l\) are homogeneous of degree \(l+1\), there exist polynomials \(p_l,q_l\in k[z]\) such that\vspace{-.32cm}
\[A_l(a,b)=a^{l+1}p_l(z)\qquad\text{and}\qquad B_l(a,b)=a^{l+1}q_l(z)\]
A direct computation gives
\[D(a)=b^2=a^2z^2\qquad\text{and}\qquad D(b)=a^2\]
The unique extension of \(D\) to the rational function field \(k(a,b)\) satisfies
\vspace{.06cm}\[
D(z)=D(\frac{b}{a})=\frac{D(b)a-bD(a)}{a^2}=\frac{a^3-b^3}{a^2}=a-bz^2=a(\smash{1-z^3})
\vspace{.06cm}\]
The recurrences for \(A_l\) and \(B_l\) give \(p_0=1\), \(q_0=z\) and
\EqL{p_{l+1}(z)=(\smash{1-z^3})p^\prime_{l}(z)+(l+1)z^2p_l(z)\\
q_{l+1}(z)=(\smash{1-z^3})q^\prime_{l}(z)+(l+1)z^2q_l(z)}{Rec0}
By symmetry under the exchange of \(a\) and \(b\) we also have
\[A_l(a,b)=B_l(b,a)\]
and hence\vspace{-.32cm}
\EqL{q_l(z)=z^{l+1}p_l(1/z)}{Reciprocal}
\subparagraph*{}
Let \(r_l,s_l\in\Set{0,1,2}\) be determined by \(r_l\equiv-l\bmod3\) and \(s_l\equiv(1-l)\bmod3\).
The congruences in \eqref{CongruenceMod3} allow us to write
\EqL{p_l(z)=z^{r_l}P_l(\smash{z^3})}{plFactorization}
By the homogeneity of \(A_l\) we have
\[\deg P_l \le \Floor{\frac{l+1-r_l}{3}}=\Floor{\frac l3}\]
The recurrence defining \(A_{l+1}=D(A_l)\) shows by induction that the monomials \(a^{s_l}b^{l+1-s_l}\) and \(a^{l+1-r_l}b^{r_l}\) corresponding to the leading and constant terms of \(P_l\) have nonzero coefficients.
Indeed, this holds for \(A_0=a\), and applying \(D\) to these monomials produces the corresponding ones at the next step.
Since all coefficients are non-negative integers, no cancellation can occur.
Therefore\vspace{-.25cm}
\[\deg P_l = \Floor{l/3}\]
\subparagraph*{}
Write \(n_l=\Floor{l/3}\) and consider the reciprocal polynomial
\[P_l^\ast(t)=t^{n_l}P_l(1/t)\]
where \(t=z^3\).
Using Equation \eqref{Reciprocal},
\Eq{q_l(z)&=z^{l+1-r_l}P_l(\smash{1/t})\\
&=z^{l+1-r_l-3n_l}P_l^\ast(\smash{t})\\
&=z^{s_l}P_l^\ast(\smash{z^3})}
where the last equality follows from \(s_l=l+1-r_l-3n_l\).

\begin{Pro}{}{EquivalentConj}
Using the above notation, Conjecture \ref{conj:OriginalConjecture} holds if and only if
\EqL{
\boxed{
\gcd(P_l,P_l^\ast)=1\qquad \text{for all }l\in\N_0
}}{ConjFormula2}
\end{Pro}

\vspace{.12cm}
\begin{Proof}
Since \(A_l\) and \(B_l\) are homogeneous of degree \(l+1\), we obtain\vspace{-.12cm}
\[
A_l=a^{s_l}b^{r_l}\,\widetilde P_l(a,b)
\qquad\qquad
B_l=a^{r_l}b^{s_l}\,\widetilde P_l^\ast(a,b)
\vspace{-.25cm}\]
where\vspace{-.2cm}
\[
\widetilde P_l(a,b)=a^{3n_l}P_l(\smash{b^3/a^3})
\qquad\qquad
\widetilde P_l^\ast(a,b)=a^{3n_l}P_l^\ast(\smash{b^3/a^3})
\]
Since the constant and leading coefficients of \(P_l\) are nonzero, the powers of \(a\) and \(b\) appearing explicitly in the above expressions of \(A_l\) and \(B_l\) are their exact multiplicities.

Let \(R_l=\gcd(P_l,P_l^\ast)\).
Then
\[
\gcd(A_l,B_l)\sim
(ab)^{\min\{r_l,s_l\}}
a^{3\deg R_l}R_l(\smash{b^3/a^3})
\]
with\vspace{-.5cm}
\[
\min\{r_l,s_l\}=
\begin{cases}
1&\text{if }l\equiv2\bmod3\\
0&\text{otherwise}
\end{cases}
\]
Therefore\vspace{-.2cm}
\[
\deg\gcd(A_l,B_l)
=
2\min\{r_l,s_l\}+3\deg R_l
\]
Since \(P_0=P_0^\ast=1\), Equation \eqref{ConjFormula} holds if and only if \(\deg R_l=0\) for every \(l\in\N_0\).
\end{Proof}

\paragraph*{}
We derive the recurrence relations satisfied by \(P_l\) and \(P_l^\ast\).
Substituting \eqref{plFactorization} into the first recurrence in \eqref{Rec0} and setting \(t=z^3\) we obtain
\vspace{.06cm}
\[P_{l+1}(t)=\begin{cases}
	(l+1)P_l(t)+3(1-t)P_l^\prime(t) & \text{if }l\equiv0\bmod3\\
	(2+(l-1)t)P_l(t)+3t(1-t)P_l^\prime(t) & \text{if }l\equiv1\bmod3\\
	(1+lt)P_l(t)+3t(1-t)P_l^\prime(t) & \text{if }l\equiv2\bmod3
\end{cases}\vspace{.06cm}\]
with initial condition \(P_0(t)=1\).
In a similar way,
\vspace{.06cm}\[P_{l+1}^\ast(t)=\begin{cases}
	(1+lt)P_l^\ast(t)+3t(1-t)(P_l^\ast)^\prime(t) & \text{if }l\equiv0\bmod3\\
	(l+1)P_l^\ast(t)+3(1-t)(P_l^\ast)^\prime(t) & \text{if }l\equiv1\bmod3\\
	(2+(l-1)t)P_l^\ast(t)+3t(1-t)(P_l^\ast)^\prime(t) & \text{if }l\equiv2\bmod3
\end{cases}\vspace{.06cm}\]
with initial condition \(P_0^\ast(t)=1\).
These recurrences show that \(P_l,P_l^\ast\in\Z[t]\). We may thus regard them as polynomials over \(\R\) and study their complex roots. The first polynomials in the sequences are collected in the following table:
\[\begin{array}{c|c|c}
	l & P_l(t) & P_l^\ast(t)\\\hline
	0 & 1 & 1\\
	1 &	1 & 1\\
	2 &	2 & 2\\
	3 &	4t+2 & 2t+4\\
	4 &	4t+20 & 20t+4\\
	5 &	80t+40 & 40t+80\\
	6 & 160t^2+520t+40 & 40t^2+520t+160
\end{array}\]

\paragraph*{}
The main result of this section is the following, which, by the previous work, directly implies Conjecture \ref{conj:OriginalConjecture}:

\begin{The}{}{MainThm}
Using the above notation, the polynomials \(P_l\) and \(P_l^\ast\) have degree \(\Floor{l/3}\) and their zeros are real, negative and simple for every \(l\in\N_0\).
Moreover, for every \(l\ge 3\), the zeros of \(P_l\) and \(P_l^\ast\) strictly interlace: they are distinct, and occur in alternating order on the real line.
Consequently,\vspace{-.12cm}
\[\gcd{(P_l,P_l^\ast)}=1 \qquad\text{over }\Q[t]\]
\end{The}

\vspace{.32cm}
\subparagraph*{}
By Theorem \ref{the:MainThm}, \(P_l\) and \(P_l^\ast\) are coprime in \(\Q[t]\).
Hence there exist \(C_l,D_l\in\Q[t]\) such that
\vspace{-.4cm}\[
    C_lP_l+D_lP_l^\ast=1
\]
For any field \(k\) of characteristic zero, the canonical embedding \(\Q\hookrightarrow k\) makes this an identity in \(k[t]\). Consequently, \(P_l\) and \(P_l^\ast\) are coprime in \(k[t]\). Now Conjecture \ref{conj:OriginalConjecture} follows by Proposition \ref{pro:EquivalentConj}.

\paragraph*{}
The rest of this section is devoted to the proof of Theorem \ref{the:MainThm}.
We first introduce the following notation. Let \(f,g\in\R[z]\) be two polynomials of degree \(n\) with positive leading coefficients.
We write \(f\prec g\) if the zeros of \(f\) and \(g\) are real, simple, and satisfy
\[
0<\alpha_1<\beta_1<\alpha_2<\beta_2<\cdots
  <\alpha_n<\beta_n<1
\]
where the \(\alpha_\i\) are the zeros of \(f\) and the \(\beta_\i\) are those of \(g\).
For \(n=0\), we use the convention \(f\prec g\) for all positive constants \(f,g\).

\begin{Le}{}{}
Suppose that \(f\prec g\). Then the following properties hold:
\begin{align}
zg^\prime+cg
  &\prec (z-1)f^\prime+cf
  &&\text{for }0<c\le1
  \label{Lemma1}\\
z(z-1)f^\prime+(z-a)f
  &\prec z(z-1)g^\prime+(z-b)g
  &&\text{for }0<a<b<1
  \label{Lemma2}
\end{align}
\end{Le}

\vspace{.2cm}
\begin{Proof}
The proof uses properties of interlacing in terms of the signs of the residues of a rational function. Since we haven't found a precise reference to this result in the literature, we write a proof in detail.
For the relation between interlacing and signs in partial fraction decompositions, see \cite[Section~1.3, Lemma~1.20]{Fisk}.
\\
If \(f\) and \(g\) are positive constants, both properties hold trivially, so assume \(n\ge1\).
Write
\[
R(z)=\frac{f(z)}{g(z)}
\]
The partial fraction decomposition of \(R\) takes the form
\[
R(z)=R(0)+\sum_{\i=1}^{n}
        \frac{c_\i z}{\beta_\i(z-\beta_\i)}
\qquad\text{with }
c_\i=\frac{f(\beta_\i)}{g^\prime(\beta_\i)}>0
\]
where the \(c_\i\) are positive because the zeros of \(f\) and \(g\) are alternating.
Indeed, writing \(f=A\*\prod_{\j=1}^{n}(z-\alpha_\j)\) and  \(g=B\*\prod_{\j=1}^{n}(z-\beta_\j)\) one has
\[g^\prime(\beta_\i)=B\prod_{\substack{1\leq \j\leq n\\ \j\neq\i}}(\beta_\i-\beta_\j)\]
and therefore\vspace{-.32cm}
\[c_\i=\frac{A}{B}\cdot \frac{\displaystyle\prod_{\j=1}^n(\beta_\i-\alpha_\j)}{\displaystyle\prod_{\substack{1\le\j\le n\\ \j\ne\i}}(\beta_\i-\beta_\j)}\]
In both products, the negative factors are those with \(\j>\i\).
Thus the number of negative factors is the same, and \(c_\i>0\).
Moreover \(R(0)>0\), since \(f(0)\) and \(g(0)\) both have sign \((-1)^n\).

For a given \(d\in\R\), \(0<d\le1\) define the rational function
\[
H_d(z)=z(z-1)R^\prime(z)+dR(z)
\]
Using the partial fraction decomposition of \(R\), a direct computation shows that
\[
H_d(z) = dR(0)+ \sum_{\i=1}^{n}
\frac{
c_\i z\bigl[(1-z)\beta_\i(1-d)+zd(1-\beta_\i)\bigr]
}{
\beta_\i(z-\beta_\i)^2
}
\]
In particular,
\EqL{H_d(z)>0\qquad \text{for all } z\in(0,1)\setminus\{\beta_1,\dots,\beta_n\}}{eq:Positivity}

To prove \eqref{Lemma1}, we take
\[
A(z)=(z-1)f^\prime(z)+cf(z)
\qquad\text{and}\qquad
B(z)=zg^\prime(z)+cg(z)\]
These polynomials have degree \(n\) and positive leading coefficient.
From the identities
\[
B(0)=cg(0)
\qquad
B(\beta_\i)=\beta_\i\,g^\prime(\beta_\i)
\]
we obtain that \(B\) has exactly one simple zero in each of the following intervals:
\[
(0,\beta_1),\;
(\beta_1,\beta_2),\dots,
(\beta_{n-1},\beta_n)
\]
Denote these zeros by \(u_1<\cdots<u_n\).
From \(B(u_\i)=0\) we obtain
\[
\frac{g^\prime(u_\i)}{g(u_\i)}=-\frac{c}{u_\i}
\]
and hence\vspace{-.2cm}
\[
\frac{A(u_\i)}{g(u_\i)} = (u_\i-1)R^\prime(u_\i)+\frac{c}{u_\i}R(u_\i) = \frac{H_c(u_\i)}{u_\i}>0
\]
It follows that the values \(A(u_\i)\) alternate in sign.
To complete the proof, we note that
\[
A(u_n)<0
\qquad\text{and}\qquad
A(1)=cf(1)>0
\]
and hence \(A\) must have a zero in each of the intervals
\[
(u_1,u_2),\dots,(u_{n-1},u_n),(u_n,1)
\]
Since \(\deg A=n\), these are all the zeros of \(A\), and in particular they are simple. By definition, it follows that \(B\prec A\).
\\
The proof of \eqref{Lemma2} is analogous, considering the polynomials
\Eq{
A(z)&=z(z-1)f^\prime(z)+(z-a)f(z)
\\
B(z)&=z(z-1)g^\prime(z)+(z-b)g(z)\qedhere
}
%
\end{Proof}

\vspace{.05cm}
\begin{Proof}[\bf of Theorem \ref{the:MainThm}]\hspace{-.5cm}
Evaluating the recurrence for \(P_l\) at \(t=1\), it is immediate to check that\vspace{-.32cm}
\[
P_l(1)=P_l^\ast(1)=l!
\]
Using the same notation \(n_l=\lfloor l/3\rfloor\) as above, define
\[F_l(z)=(z-1)^{n_l} P_l\!\left(\frac{z}{z-1}\right)
\qquad\text{and}\qquad
G_l(z)=(z-1)^{n_l}P_l^\ast\!\left(\frac{z}{z-1}\right)
\]
Multiplication by \((z-1)^{n_l}\) clears the denominators, so \(F_l,G_l\in\R[z]\).
Their leading coefficients are \(P_l(1)=P_l^\ast(1)=l!>0\), and hence both have degree \(n_l\).
Applying the chain rule to the recurrences for \(P_l\) and \(P_l^\ast\) we obtain
\Eq{F_{l+1}=\begin{cases}
	F_l+3(z-1)F_l^\prime &\text{if }l\equiv0\bmod3\\
	2F_l+3zF_l^\prime &\text{if }l\equiv1\bmod3\\
	(3z-1)F_l + 3z(z-1)F_l^\prime &\text{if }l\equiv2\bmod3
\end{cases}
\\[-5pt]
\tag{\(\ast\)}
\\[-5pt]
G_{l+1}=\begin{cases}
	G_l+3zG_l^\prime&\text{if }l\equiv0\bmod3\\
	2G_l+3(z-1)G_l^\prime &\text{if }l\equiv1\bmod3\\
	(3z-2)G_l+3z(z-1)G_l^\prime &\text{if }l\equiv2\bmod3
\end{cases}
}

We prove by induction that
\EqL{
\begin{cases}
F_l\prec G_l &\text{if } l\equiv 0,2\bmod 3\\
G_l\prec F_l &\text{if } l\equiv 1\bmod 3
\end{cases}
}{Alternating}
For \(l=0\) there is nothing to prove, since \(F_0=G_0=1\).
Assume that the assertion holds for some \(l\ge0\). If \(l\equiv 0\bmod 3\), the inductive hypothesis gives \(F_l\prec G_l\). Applying \eqref{Lemma1} with \(c=1/3\) and multiplying both polynomials by 3 we immediately obtain
\[
G_{l+1}\prec F_{l+1}
\]
If \(l\equiv 1\bmod 3\) we have \(G_l\prec F_l\). Applying \eqref{Lemma1} with \(c=2/3\), \(f=G_l\) and \(g=F_l\) we immediately obtain\vspace{-.2cm}
\[
F_{l+1}\prec G_{l+1}
\]
Finally, if \(l\equiv 2\bmod 3\) we have $F_l\prec G_l$.
Applying \eqref{Lemma2} with \(a=1/3\) and \(b=2/3\) we obtain again\vspace{-.2cm}
\[
F_{l+1}\prec G_{l+1}
\]

For \(l\ge3\), this shows that the zeros of \(F_l\) and \(G_l\) are simple, belong to the open interval \((0,1)\) and strictly interlace.
To complete the proof, note that
\[\NMapH{(0,1)}{(-\infty,0)}{z}{\dfrac{z}{z-1}}{4pt}\]
is a decreasing bijection.
Under this map, the zeros of \(F_l\) and \(G_l\) are transported onto the zeros of \(P_l\) and \(P_l^\ast\).
In particular, this also shows that the zeros of the latter polynomials are simple and strictly interlace.
Since \(P_l\) and \(P_l^\ast\) have no common zero in \(\C\), we obtain
\[\gcd(P_l,P_l^\ast)=1 \qquad\quad
       \text{over }\Q[t]\text{ for all }l\in\N_0\qedhere\]
\end{Proof}


\section{Thomas decomposition and computational verification}\label{Thomas}
When studying a system of polynomial equations and inequations, the degrees of the polynomials alone are generally not sufficient to understand the structure of its solution set: specializing some variables may lower the degree in another variable or cause distinct roots to coincide.
For example, the product of the degrees of two projective plane curves without a common component gives their total intersection multiplicity \cite[Chapter I, Section 1.6]{Shafarevich} but only an upper bound on the number of distinct affine intersection points, since some intersections may lie at infinity or have multiplicity greater than 1.
These phenomena are particularly relevant to geometric degree computations, which involve counting the distinct solutions of zero-dimensional systems.

An \textbf{(algebraic) Thomas decomposition} splits a polynomial system \(S\) into a finite collection of systems \(S_1,\dots,S_r\) of equations and inequations such that
\[
Z(S)=\bigsqcup_{\i=1}^r Z(S_\i)
\]
where the \(S_\i\) satisfy additional conditions that allow their distinct solutions to be counted from the degrees of their equations when \(Z(S)\) is finite.

\paragraph*{}
Work over the polynomial ring \(k[u_1,\dots,u_n]\) with the fixed ordering \(u_1<\dots<u_n\) on the indeterminates. Let \(p\) be a nonconstant polynomial. The \textbf{leader} of \(p\) is the largest variable appearing in \(p\), and is denoted by \(\ld{p}\).
If \(\ld{p}=u_\i\) we can write
\[p=c_d(u_1,\dots,u_{\i-1})u_\i^d+\dots+c_0(u_1,\dots,u_{\i-1})\in k[u_1,\dots,u_{\i-1}][u_\i]\]
with \(c_d\ne0\).
The integer \(d\ge1\) and the polynomial \(\init{p}=c_d\) are respectively called the \textbf{degree in the leader} and the \textbf{initial} of \(p\).
By convention, we set \(\ld{p}=1\) if \(p\) is constant.
When considering one equation \(p=0\), we speak of the leader, the degree and the initial of the equation, instead of the leader, the degree and the initial of \(p\); the same applies when considering one inequation \(q\ne 0\).
\subparagraph*{}
A \textbf{triangular chain} is a finite sequence of equations with pairwise distinct leaders.
We will always assume that the equations in a triangular chain are nontrivial (the polynomials defining them are nonconstant) and ordered by increasing leaders.
It follows immediately from the definition that a triangular chain in \(k[u_1,\dots,u_n]\) consists of at most \(n\) equations.
The indeterminates that are not the leader of any equation in a triangular chain are called \textbf{parametric}.

\begin{Exa}{}{}
The sequence \((x=0, yz=0)\) is a triangular chain in \(k[x,y,z]\) with \(x<y<z\).
The variable \(y\) is parametric for such a chain.
\end{Exa}

\vspace{.4cm}
\noindent
Consider the chain in the preceding example and the projection
\[\MapH{\pi}{Z(x,yz)}{\A^2}{(a,b,c)}{(a,b)}{0pt}\]
Then \(\pi^{-1}(a,b)\) consists of all points in \(Z(x,yz)\) such that the first coordinates are equal to \(a,b\).
The affine variety \(Z(x,yz)\) can be decomposed using the fibers of \(\pi\):\vspace{.2cm}
\[Z(x,yz)=\bigsqcup_{(a,b)\in \A^2}\pi^{-1}(a,b)=\Set{(0,0,c)\colon c\in K}\;\;\sqcup\bigcup_{b\in\A^1\setminus\Set0}\Set{(0,b,0)}
\]
The first set is the solution set of \(x=0\), \(y=0\) and the second is the solution set of \(x=0\), \(z=0\) and \(y\ne0\). These systems are obtained using vanishing and non-vanishing conditions on the initials of the equations in the original system. If one initial vanishes, the degree of the corresponding equation drops.
We give two further examples, which illustrate this splitting for a finite solution set and for a curve with a parametric variable:

\begin{Exa}{}{Thomas1}
Work over \(k[u,v]\) with \(u<v\), and consider the system
\Eq{S\colon\begin{cases}f=u(\smash{u^2-1})=0\\
g=uv^2+v=0\end{cases}}
The first equation gives \(u\in\Set{-1,0,1}\).
The second equation has initial \(u\) and leader \(v\), so its degree in \(v\) drops when the lower variable \(u\) is specialized to \(0\).
Proceeding as before, consider the projection
\[\MapH{\pi}{Z(S)}{\A^1}{(a,b)}{a}{0pt}\vspace{.12cm}\]
and write \(F_a=\pi^{-1}(a)\) for each \(a\in\A^1\).
Then\vspace{.06cm}
\EqL{Z(S)
&=F_{0}\sqcup (F_{-1}\cup F_{1})\\[2pt]
&=\Set{(0,0)} \,\sqcup\, \Set{(-1,0), (-1,1), (1,0), (1,-1)}}{Disj1}
The first set is the solution set of \(u=0\), \(v=0\) and the second is the solution set of \(u^2-1=0\), \(g=0\).
Note that no inequation on the initial is used in the second system, since the equation \(u^2-1=0\) already forces \(u\ne 0\).
The fiber over \(u=0\) has cardinality 1, while the fibers over \(u=-1\) and \(u=1\) both have cardinality 2.
The two disjoint subsets in \eqref{Disj1} have constant fiber cardinalities over their respective images under \(\pi\).
\end{Exa}

\begin{Exa}{}{Thomas2}
We slightly modify the preceding example. Consider the affine plane curve of equation
\[g=uv^2+v=0\]
The indeterminate \(u\) is now parametric.
Solving the equation for \(v\), we obtain
\EqL{Z(g)&=\Set{(0,0)}\sqcup \Set{\smash{(a,0), (a,-1/a)\mid a\in \A^1\setminus\Set{0}}}}{Disj2}
The first set is the solution set of \(u=0\) and \(v=0\); the second is the solution set of \(g=0\) and \(u\ne0\).
Under the projection onto the \(u\)-coordinate, the fiber over \(u=0\) has cardinality 1, while the fibers over \(u\ne 0\) have cardinality 2.
\end{Exa}

\paragraph*{}
It is convenient to present an algebraic system \(S\) as an ordered pair \((E,I)\) of two vectors consisting of the polynomials defining its equations and inequations, respectively.
We write \(f\in S\) to express that \(f\in E\cup I\).
\\
For a variable \(v\), let \(S_{\le v}\) denote the subsystem consisting of the equations and inequations in \(S\) whose leaders are less than or equal to \(v\).
Similarly, \(S_{<v}\) denotes the system obtained by only considering the equations and inequations in \(S\) whose leaders are strictly less than \(v\).
We regard \(S_{\le u_\i}\) and \(S_{<u_\i}\) as systems in \(k[u_1,\dots,u_\i]\) and \(k[u_1,\dots,u_{\i-1}]\) respectively.
The ordering on the variables is extended by setting \(1<u_1\). Using this convention, \(S_{\le1}\) consists of the constant equations and inequations in \(S\); in particular, it is the empty system whenever \(S\) contains no constant polynomials.
\subparagraph*{}
We say that \(S\) is \textbf{quasi-simple} if the following conditions hold:
\begin{enumerate}[leftmargin=0.8cm, topsep=0.28cm]
	\item The polynomials in \(S\) are nonconstant with pairwise distinct leaders.
	\item The equation \(\init{h}=0\) has no common solution with \(S_{\le\ld{\init{h}}}\) for every \(h\in S\).
\end{enumerate}

\noindent
The second condition is imposed on the solutions of the subsystem determined by the leader of the initial, not on those of the full system. Indeed, solving a triangular system, we proceed by successively specializing the lower variables and solving the corresponding equations in their leaders.

The following proposition explains the relation between the two formulations in the case of a quasi-simple system:

\begin{Pro}{}{Extending}
Let \(S=(E,I)\) be a quasi-simple system in \(k[x_1,\dots,x_n]\).
The system \(S\) has a solution, and for every \(l\in\{1,\dots,n\}\), any solution of \(S_{\le x_l}\) extends to a solution of \(S\).
Consequently, for every polynomial \(g\in k[x_1,\dots,x_n]\), the following conditions are equivalent:
\begin{enumerate}[leftmargin=0.8cm, topsep=0.28cm]
    \item The equation \(g=0\) has no common solution with \(S_{\le\ld{g}}\).
    \item The equation \(g=0\) has no common solution with \(S\).
\end{enumerate}
\end{Pro}

\begin{Proof}
We first prove the extension property: every element \((a_1,\dots,a_l)\in Z(S_{\le x_l})\) can be extended successively to the variables \(x_{l+1},\dots,x_n\).
Since the leaders of the polynomials in \(S\) are pairwise distinct, there is at most one \(f\in S\) with leader \(x_{l+1}\).
If no such polynomial exists, we may choose \(a_{l+1}\in K\) arbitrarily.
Otherwise, consider \(f\) as a univariate polynomial in \(x_{l+1}\). By quasi-simplicity, the equation \(\init{f}=0\) has no common solution with \(S_{\le\ld{\init{f}}}\).
Since the variable \(x_{l+1}\) does not appear in \(\init{f}\) and since \((a_1,\dots,a_{l})\) is a solution of \(S_{\le x_{l}}\) we have
\[
\init{f}(a_1,\dots,a_{l})\neq 0
\]
Consequently, the specialization \(f(a_1,\dots,a_{l},x_{l+1})\) is a nonconstant univariate polynomial of the same degree in \(x_{l+1}\) as \(f\), and hence has only finitely many roots.
If \(f\in E\), this polynomial has at least one root in the algebraically closed field \(K\); we choose \(a_{l+1}\) to be one of these roots.
If \(f\in I\), we choose \(a_{l+1}\) outside its zero set, which is possible because \(K\) is infinite.
In either case, the resulting tuple \((a_1,\dots,a_{l+1})\) satisfies \(S_{\le x_{l+1}}\).
Proceeding inductively, we obtain a solution of the full system \(S\) extending the original partial solution.
Applying the same argument from \(x_1\) onward and starting from the empty tuple also shows that \(S\) has a solution.

For the equivalence of the two conditions, note that the implication from 2 to 1 follows from the extension property.
The converse follows by restricting any common solution of \(S\) and \(g=0\) to the variables not greater than \(\ld{g}\).
\end{Proof}

\begin{Rmk}{}{RmkEx}
Non-vanishing of the initials on the solutions of the full system cannot replace the lower subsystem conditions in the definition of quasi-simplicity. These lower subsystem conditions are also needed to guarantee the extension property.
For example, consider the system\vspace{-.12cm}
\[
S=(\smash{x^2-1=0,\, xy^2-y=0,\, yz-1=0})
\vspace{-.12cm}\]
in \(\Q[x,y,z]\) with \(x<y<z\).
The initial of the third equation is \(y\), and the  corresponding lower subsystem is\vspace{-.4cm}
\[
S_{\le y}= (\smash{x^2-1=0,\ xy^2-y=0})
\]
The equation \(y=0\) has the common solution \((1,0)\) with \(S_{\le y}\) and \(S\) is not quasi-simple.
It is immediate to check that \((1,0)\) cannot be extended to a solution of the full system, thus the extension property fails.
\end{Rmk}

\vspace{-.12cm}
\paragraph*{}
Quasi-simplicity alone is not sufficient if we are interested in counting solutions without multiplicity: the system \(S=(\smash{v^2-u}=0)\) in \(\Q[u,v]\) with \(u<v\) is quasi-simple, but setting \(u=0\), the equation becomes \(v^2=0\), which has only one solution in \(\overline\Q\) with multiplicity 2.
Thus, after specialization, it may happen that some equation has multiple roots, yielding fewer geometric points than its degree.
To account for this problem, we also need to consider a square-freeness condition on the system. This condition can be expressed in terms of the discriminants of the equations and inequations.
\subparagraph*{}
Let \(h\in k[u_1,\dots,u_n]\) be a polynomial, with \(\ld{h}=u_\i\) and \(\deg_{u_\i}h=d>1\).
Its \textbf{discriminant} with respect to \(u_\i\) is the polynomial
\vspace{.06cm}
\[
\disc{h}=(-1)^{d(d-1)/2}\frac{\res{h,\sep{h}}{u_\i}}{\init{h}}\in k[u_1,\dots,u_{\i-1}]
\vspace{.04cm}\]
Here \(\sep{h}=\partial h/\partial u_\i\) is the \textbf{separant} of \(h\), and \(\res{h,\sep{h}}{u_\i}\) denotes the \textbf{resultant} of \(h\) and \(\sep{h}\) with respect to \(u_\i\).
We set by convention \(\disc{p}=1\) if \(p\) has degree at most 1 in its leader.
In particular, provided that the initial of \(h\) does not vanish, the specialization of \(h\) has \(d\) distinct roots in \(K\) if and only if its discriminant does not vanish \cite[Section 3.6, Corollary 4 and Proposition 6]{CLO}.
\subparagraph*{}
A quasi-simple system \(S=(E,I)\) is called \textbf{simple} if, for every \(h\in S\), the equation \(\disc{h}=0\) has no common solution with \(S_{\le\ld{\disc{h}}}\).

\begin{Exa}{}{}
The system \eqref{EqP} is simple with respect to any ordering in which \(x_1\) and \(x_2\) are lower than all higher jet variables: the equations are monic and linear in their pairwise distinct leaders, so all initials and discriminants are equal to 1.
\end{Exa}

\begin{Rmk}{}{}
The non-vanishing of initials ensures that the degree of every equation and inequation is preserved under every specialization of the lower variables satisfying the relevant lower subsystem.
Such specializations are called \textbf{admissible}.
On the other hand, the non-vanishing of the discriminants ensures that the specialized polynomials have distinct roots.
In particular, an inequation of degree \(d\) in its leader excludes exactly \(d\) values in every admissible fiber.
\subparagraph*{}
Once quasi-simplicity is certified, the discriminant condition can be checked on either the solution set of the relevant lower system or of the full system.
The extension property in Proposition \ref{pro:Extending} does not require characteristic zero.
The discriminant condition for the absence of multiple roots also holds in arbitrary characteristic, provided that the degree is preserved under specialization.
\end{Rmk}

\paragraph*{}
Returning to the example \(S=(\smash{v^2-u})\), we have \(\disc{\smash{v^2-u}}=4u\).
The equation \(u=0\) has solutions in the relevant lower subsystem \(S_{\le u}\) (which is the empty system, satisfied by every element of \(K\)): \(S\) is quasi-simple, but not simple.
We then proceed considering the following splitting:
\[
S_1=(u=0,v=0)\qquad
S_2=(\smash{v^2-u=0, u\ne 0})
\]
Both \(S_1\) and \(S_2\) are simple systems.
The first system describes the exceptional fiber over \(u=0\) consisting of one point, while the second one describes the fibers over \(u\neq 0\), each consisting of two distinct points.
This is an example of a Thomas decomposition of \(S\), and motivates the following general definition.
\subparagraph*{}
Let \(S=(E,I)\) be a system of polynomial equations and inequations in \(k[u_1,\dots,u_n]\). A \textbf{Thomas decomposition} for \(S\) is a finite (possibly empty) collection \(S_1,\dots,S_N\) of simple systems such that
\[Z(S)=\bigsqcup_{\i=1}^N Z(S_\i)\]

\begin{Rmk}{}{}
A Thomas decomposition of an algebraic system can be obtained algorithmically \cite[Section 2.2.1, Theorem 2.2.30]{Robertz2014}. An implementation of the algorithm is available in Maple through the package \texttt{AlgebraicThomas}. An experimental implementation using OSCAR is also available \cite{MyGitHub}.
\end{Rmk}

The triangular structure of a simple system, together with its non-vanishing conditions, allows us to study its solutions iteratively, controlling the behavior of the fibers over the lower variables.

\subsection{Zero-dimensional algebraic systems}
An algebraic system \(S\) is called \textbf{zero-dimensional} if its solution set \(Z(S)\subseteq K^n\) is finite (possibly empty).
Let \(S_1,\dots,S_N\) be a Thomas decomposition of a zero-dimensional algebraic system \(S\).
Each \(S_\i\) is in turn zero-dimensional, since \(Z(S_\i)\subseteq Z(S)\).
We can recognize zero-dimensionality using Thomas decomposition:

\begin{Le}{}{QS-ZD}
A quasi-simple system \(T\) is zero-dimensional if and only if each variable \(u_\j\) is the leader of one equation in \(T\).
If that is the case, \(T\) contains no inequations.
\end{Le}

\begin{Proof}
We first suppose that every variable is the leader of one equation in \(T\).
By quasi-simplicity, we have exactly one equation for each variable and no inequations.
These equations can be solved successively, and this gives only finitely many solutions.
Thus \(T\) is zero-dimensional.
\\
Conversely, suppose that a variable \(u_\j\) is not the leader of one equation in \(T\); let \(\j\) be the smallest index with this property.
The lower subsystem \(T_{<u_\j}\) is quasi-simple and has a solution. Fix such a solution, taking the empty tuple if \(\j=1\).
If \(u_\j\) is not the leader of any polynomial in \(T\), then it can be chosen arbitrarily.
If \(u_\j\) is the leader of one of the inequations in \(T\) then, by quasi-simplicity, this inequation is unique and its initial does not vanish after specializing its coefficients (as a univariate polynomial in \(u_\j\)).
We thus obtain a nonzero univariate polynomial, and only finitely many possibilities for \(u_\j\) are excluded.
In both cases, since an algebraically closed field of characteristic zero is infinite, \(u_\j\) can take infinitely many values.
By Proposition \ref{pro:Extending}, each of these choices produces a distinct solution of \(T\).
Thus \(Z(T)\) is infinite, and \(T\) is not zero-dimensional.
\end{Proof}

\begin{Pro}{}{}
An algebraic system \(S\) is zero-dimensional if and only if every simple system in a Thomas decomposition of \(S\) contains only equations, one for each possible variable.
\end{Pro}

\vspace{.12cm}
\begin{Proof}
Since the solution sets of the simple systems in a Thomas decomposition form a finite partition of \(Z(S)\), the system \(S\) is zero-dimensional if and only if each of these systems has the same property.
The claim follows from the previous lemma.
\end{Proof}

\vspace{.012cm}
\begin{Rmk}{}{RmkEx2}
Both implications of the preceding lemma may fail without the assumption of quasi-simplicity.
Work in \(\Q[x,y,z]\) with \(x<y<z\).
The system\pagebreak\noindent
\[S\colon\begin{cases}
	x-1=0\\
	y^2-x=0\\
	(x-1)z^2+(y-1)z+(x-1)=0
\end{cases}
\]
is not quasi-simple: the initial \(x-1\) of the third equation has a common solution with the relevant lower subsystem.
It is immediate to check that \(Z(S)\) contains the line \(\Set{(1,1,z)\colon z\in K}\).
Thus \(S\) is not zero-dimensional, although every variable is the leader of one equation in \(S\).

Conversely, consider the system
\[S^\prime\colon\begin{cases}
	y^2-x=0\\
	y-x^2=0
\end{cases}\]
in \(\Q[x,y]\) with \(x<y\).
Both equations have leader \(y\), so this system is not quasi-simple.
Substituting \(y=x^2\) into the first equation gives \(x^4-x=0\), which has four distinct roots in \(K\). Thus \(S^\prime\) is zero-dimensional, although no equation has leader \(x\).
\end{Rmk}

\paragraph*{}
Fix a simple system \(S_\i\) in the given Thomas decomposition of \(S\).
Denote by \(f_{\i\j}=0\) the equation in \(S_\i\) with leader \(u_\j\) and write \(d_{\i\j}=\deg_{u_\j}f_{\i\j}\).
Consider the projections
\[\MapH{\pi_\j}{Z(f_{\i1},\dots,f_{\i\j})}{Z(f_{\i1},\dots,f_{\i(\j-1)})}{(a_1,\dots,a_\j)}{(a_1,\dots,a_{\j-1})}{0pt}
\]
Since \(S_\i\) is simple and there is no inequation by the preceding proposition, every fiber contains exactly \(d_{\i\j}\) points.
Thus by induction
\[|Z(S_\imath)\!|=\prod_{\j=1}^{n}d_{\i\j}\]
Since the solution sets of \(S_1,\dots,S_N\) give a partition of the solution set of \(S\), we obtain the formula\vspace{-.2cm}
\EqL{|Z(S)\!|=\sum_{\i=1}^N\prod_{\j=1}^{n}d_{\i\j}}{FormulaCountingPoints}

\pagebreak\noindent
\begin{Rmk}{}{}
The original zero-dimensional system \(S\) may contain inequations.
Thus the product of the degrees of the equations may not directly give the number of solutions before the system has been decomposed into simple systems.
For example,
\[
\begin{cases}
	u^2-1=0\\
	u-1\neq 0
\end{cases}
\]
has only the solution \(u=-1\), although its equation has degree 2.
The Thomas decomposition algorithm produces the simple system consisting only of the equation \(u+1=0\), for which the product formula gives the correct number of solutions.
\end{Rmk}

\paragraph*{}
The companion notebook uses Formula \eqref{FormulaCountingPoints} to count the distinct affine intersection points in the experiments described in Section \ref{Motivation}.

\subsection{Certifying hyperplane sections}
We return to the system \(\Sigma_2\) and use the notation of Subsection \ref{Infinity}.
For a fixed \(m\in\N_0\) consider
\vspace{-.4cm}
\Eq{P&=A_{m+1}+U\\
Q&=B_{m+1}+V\qquad\quad \text{for some }U,V\in W_m}
and set \(d=m+2\), \(H=\gcd(A_{m+1},B_{m+1})\) and \(r=\deg H\).
Assume that \(P,Q\) are coprime, that \(H\) is square-free, and that\vspace{-.2cm}
\[
\gcd(H,J_a,J_b)=1
\]
The projective curves of equations \(P^h=0\) and \(Q^h=0\) then have \(r=\deg H\) intersection points at infinity, each with intersection multiplicity 1.
By Theorem \ref{the:TheConj} and Bézout's theorem,\vspace{-.2cm}
\[|Z(P,Q)\!| \le \sum_{p\in Z(P,Q)} I_p(\smash{P^h,Q^h}) =d^2-r =\deg\X^{(m)}
\]
where \(I_p(\smash{P^h,Q^h})\) denotes the local intersection multiplicity at \(p\), with equality if and only if every affine intersection point has multiplicity 1.

\paragraph*{}
Applying Formula \eqref{FormulaCountingPoints} to a Thomas decomposition of the system \(P=0\), \(Q=0\) gives the exact number \(|Z(P,Q)\!|\) of distinct affine intersection points.
If this quantity coincides with \(d^2-r\), we certify that every affine intersection is simple.
Together with the conditions already verified at infinity, this shows that the chosen pair of hyperplane sections realizes the geometric degree.

\begin{Rmk}{}{NotSQFree}
As will be shown in Section \ref{Further}, this criterion may not apply to more general systems even for generic hyperplane sections.
\end{Rmk}


\section{Further applications and extensions}\label{Further}
The approach presented in Section \ref{Motivation} can be applied to more general systems of the form\vspace{-.2cm}
\[\begin{cases}
	x_1^\prime = \varphi(x_1,x_2)\\
	x_2^\prime=\psi(x_1,x_2)
\end{cases}\]
Writing \(x_1=a\) and \(x_2=b\), consider the derivation
\[D_{\varphi,\psi}=\varphi(a,b)\frac{\partial}{\partial a}+\psi(a,b)\frac{\partial}{\partial b}\]
on \(k[a,b]\). For this system, define \(A_0=a\), \(B_0=b\) and
\Eq{A_{l+1}&=D_{\varphi,\psi}(A_l)\\
B_{l+1}&=D_{\varphi,\psi}(B_l)\qquad(l\ge 0)}
The equations of the \(m\)-th prolongation are then equivalent to \(x_1^{(l)}=A_l(x_1,x_2)\) and \(x_2^{(l)}=B_l(x_1,x_2)\) for \(l\in\Set{1,\dots,m+1}\).

If \(\varphi,\psi\) are homogeneous of the same degree \(s\ge 2\), each nonzero \(A_l,B_l\) is homogeneous of degree \(1+l(s-1)\).
If the two polynomials of order \(m+1\) are coprime, the corresponding projective curves have no intersection points at infinity, and the associated prolongation variety has degree \(d^2\), where
\[d=1+(m+1)(s-1)\]
When the GCD is nonconstant, the argument of Theorem \ref{the:TheConj} can still be applied if its hypotheses are satisfied.
Otherwise, the intersection multiplicities at infinity must be studied more carefully.

\paragraph*{}
We now specialize to the prolongation varieties associated with the systems \(\Sigma_n\) (\(n\ge2\)).
Write
\[W_m=\text{Span}_k(1,A_{n,0},B_{n,0},\dots,A_{n,m},B_{n,m})\qquad\quad\text{for each }m\in\N_0\]
By symmetry under exchanging \(a\) and \(b\), it follows that \(B_{n,l}(a,b)=A_{n,l}(b,a)\) for every \(l\ge0\).

As in Subsection \ref{Infinity}, we consider normalized hyperplane sections
\Eq{
P&=A_{n,m+1}+U\\
Q&=B_{n,m+1}+V
\qquad (U,V\in W_m)
}
For \(l=1,2,3\) we have
\Eq{A_{n,1}&=b^n\\
A_{n,2}&=na^{n}b^{n-1}\\
A_{n,3}&=na^{n-1}b^{n-2}(\smash{(n-1)a^{n+1}+nb^{n+1}})}
For a monomial \(a^sb^t\), the derivation
\[
D_n=b^n\frac{\partial}{\partial a}+a^n\frac{\partial}{\partial b}
\]
satisfies
\[
D_n(\smash{a^sb^t})=sa^{s-1}b^{t+n}+ta^{s+n}b^{t-1}
\]
Every nonzero term changes both exponents by \(-1\) modulo \(n+1\).
Starting from \(A_{n,0}=a\) and \(B_{n,0}=b\), it follows by induction that every monomial \(a^sb^t\) appearing in \(A_{n,l}\) satisfies
\[
s\equiv1-l\pmod{n+1}\qquad\text{and}\qquad
t\equiv-l\pmod{n+1}
\]
while every monomial appearing in \(B_{n,l}\) satisfies
\[s\equiv-l\pmod{n+1}\qquad\text{and}\qquad t\equiv1-l\pmod{n+1}\]
\subparagraph*{}
Fix \(m\in\N_0\) and let \(r\in\Set{0,\dots,n}\) be such that \(m\equiv r\bmod{(n+1)}\). Suppose first that \(1\leq r\leq n-1\). The congruences above then show that
\[
(ab)^{n-r}\mid \gcd(A_{n,m+1},B_{n,m+1})
\]
Since both \(A_{n,m+1}\) and \(B_{n,m+1}\) are divisible by \(ab\), their common zeros at infinity include \([0:1:0]\) and \([0:0:1]\).

Now consider the intersection multiplicities at these two points.
The above congruences and the recurrence \(A_{n,l+1}=D_n(A_{n,l})\) show by induction that the least exponent of \(b\) occurring in \(A_{n,l}(1,b)\) is the representative of \(-l\bmod{(n+1)}\) in \(\Set{0,\dots,n}\). For \(l=m+1\), this exponent is \(n-r\).
Moreover, the polynomial \(B_{n,m+1-r}\in W_m\) contributes a term proportional to \(Z^{(n-1)r}\) after homogenization in the chart \(a\ne0\).
These terms suggest an intersection multiplicity of \((n-1)r(n-r)\) at \([0:1:0]\) for a generic choice of the coefficients.
By symmetry, the same multiplicity is expected at \([0:0:1]\).
For \(r=0\) or \(r=n\), these two points are not common zeros of the homogenized polynomials.

We propose the following generalization of the formulas for the geometric degree in~\cite{Pogudin2025}:

\begin{Conj}{}{}
For every \(n\ge2\) and \(m\ge0\), let \(r\in\Set{0,1,\dots,n}\) be such that \(m\equiv r\bmod{(n+1)}\).
Then\vspace{-.24cm}
\[\gcd{(A_{n,m+1},B_{n,m+1})}\sim \begin{cases}
	(ab)^{n-r} & \text{if }1\leq r\leq n-1 \\
	1 & \text{if }r\in\Set{0,n}
\end{cases}\]
and\vspace{-.2cm}
	\[\deg\X_n^{(m)}=((n-1)m+n)^2-2(n-1)r(n-r)\]
\end{Conj}

\vspace{.4cm}\noindent
The formula for the degree would follow from the GCD statement, provided that the expected local intersection multiplicities are attained generically.

For \(n=3\), the cases \(m=1,2\) illustrate two different reasons why Theorem~\ref{the:TheConj} does not apply.
For \(m=1\), the polynomial
\[
H=\gcd(A_{3,2},B_{3,2})\sim a^2b^2
\]
is not square-free.
For \(m=2\), instead, the polynomial
\[
H=\gcd(A_{3,3},B_{3,3})\sim ab
\]
is square-free.
However, since every element of \(W_2\) has homogeneous components only in degrees \(0,1,3,5\), we have \(F_6=G_6=0\) for every choice of \(U,V\).
Consequently, \(J_a=J_b=0\), and the transversality condition of the theorem cannot be satisfied.

\paragraph*{}
Algebraic Thomas decomposition remains applicable to the corresponding polynomial systems.
For finite affine intersections, resultants offer an alternative: after one variable has been eliminated, the distinct intersection points can be counted by examining the fibers over the roots of the resulting univariate polynomial.
\subparagraph*{}
When the conditions of Theorem \ref{the:TheConj} are not satisfied, we restrict our computational verification to the conjectured GCD formula.
We verified this formula in OSCAR for every \(n\in\Set{2,\dots,100}\) and \(m\in\Set{0,\dots,1000}\).
The computations can be reproduced using the companion notebook.


\vspace{-.25cm}
\section*{Acknowledgments}\vspace{-.25cm}
This is a contribution to Project-ID 286237555 - TRR 195 - by the German Research Foundation (DFG).

I am grateful to Leonardo Lanciano for writing to inform me of the publication of his paper \cite{Lanciano2026}.
Reading it introduced me to the conjecture and led me to consider Thomas decomposition as a possible approach to address the problem.
\\
I also thank Leonardo, Gleb Pogudin and my PhD supervisors Annette Bachmayr and Daniel Robertz for their helpful comments and suggestions on the organization of this second version of the paper.

\vspace{-.4cm}
\fancypagestyle{fancyBibliography}{
	\renewcommand{\headrulewidth}{0.25pt}
	\renewcommand{\footrulewidth}{0.25pt}
	\fancyfoot[LE,RO]{}
	\fancyhead[LO,RE]{}
}
\pagestyle{fancyBibliography}
\renewcommand{\bibname}{}

\titleformat{\chapter}[display]{\bfseries\Large}	{\filleft\MakeUppercase{\chaptertitlename} \HUGE\thechapter}{.4ex}{\vspace{1ex}\filleft}[\vspace{3.5ex}]
\titlespacing*{\chapter}{0pt}{0.1\baselineskip}{0.5\baselineskip}

{
    \let\cleardoublepage\relax
	
}
\end{document}